# Dyck Numbers, III. Enumeration and bijection with symmetric Dyck paths

Gennady Eremin

***Abstract.*** Dyck paths (also balanced brackets and Dyck words) are among the most heavily studied Catalan families. This paper is a continuation of [2, 3]. In the paper we enumerate the terms of the OEIS A036991, Dyck numbers, and construct a concomitant bijection with symmetric Dyck paths. In the case of binary coding of Dyck paths we work with compact natural numbers after removing leading zeros. Analysis of binary suffixes, allowed us to obtain a bijection between arbitrary A036991 terms and symmetric A036991 terms which encode symmetric Dyck paths. The bijection generates a forest of unary non-intersecting infinite trees. The root of each bijection tree is an asymmetric term; the other nodes are symmetrical. There are an infinite number of such trees. The reader is offered a software package for working with bijection trees.

Dyck paths (balanced brackets or Dyck words) are among the most heavily studied Catalan families [1]. This article continues [2, 3]. In this paper, we enumerate the Dyck numbers, the terms of the OEIS sequences A036991 and A350346 [4], which compactly encode balanced brackets. For most of this article, we will work with the A350346 terms, a binary expansion of the A036991 terms. In our case, the data in binary code is more convenient and visual, and the processing of such information is significantly simplified.

Dyck paths are drawn in the first quadrant of the grid and start at the origin. This is natural, since the initial sections (prefixes) of such pictures are usually analyzed, and access to the initial sections is simple and convenient. But in the case of balanced brackets, we don't have a grid, we work with Dyck words. After binary coding and removal of *leading* (*left*) zeros, we get natural numbers, in which, as you know, the digits are numbered from the end, so the ending code fragments, suffixes, are more convenient and accessible. It is the analysis of a particular suffix, and often a group of suffixes, that allows you to obtain the required information about a specific Dyck number.

## 1  Levels in the OEIS A036991

In A036991 and A350346 (A036991/A350346 for short), the compactness of the data, as well as the convenience and clarity of information processing, is achieved by the natural coding of balanced brackets: the left parenthesis is replaced by zero, 0, and the right closing parenthesis is encoded by one, 1. And further, since we are working with integers (either binary or decimal), it is logical and natural to remove leading 0's (i.e. codes of the initial left parentheses, and at least one such parenthesis always exists).

In balanced brackets, the number of left and right parentheses is always the same, but when encoded, we preserve all 1's and internal 0's in binary expansion; outer (leading) zeros are removed. In this case, the balance is tracked in the last (right) fragments: *in any binary suffix, the number of 1's is not less than 0's*. If necessary, the deleted brackets are easily restored based on the number of 1's in the binary expansion. Here is a small example of balanced parenthesis encoding: $(((()))() \rightarrow 0000111011_2 = 111011_2 = 59$. As you can see, the bracket set of size (length) 2×5 = 10 is converted into a compact binary number of size 6.

A036991/A350346 terms are sorted in ascending order. Here are the initial terms (integers with different binary lengths are separated by a double slash):

A350346: 1//11//101,111//1011,1101,1111//10011,10101,10111,11011,11101,11111// 100111…
A036991: 1// 3// 5,7//11,13,15//19,21,23,27,29,31// 39,43,45,47,51,53,55,59,61,63//71,75,77,79…

Let's denote by **ĐN** the set of all Dyck numbers, the A036991/A350346 terms (regardless of the coding system, binary or decimal).

**Definition 1.** The *level* of **ĐN** is the set of Dyck numbers with the same binary length.

It is logical to associate the level number with the length of the binary code of the included terms. In our case, the first level $Ł_1 = \{1\}$, *1-level*, and the second level $Ł_2 = \{11\}$, *2-level*, contain one binary term each (decimal 1 and 3, respectively). Obviously, $\#Ł_1 = \#Ł_2 = 1$. At *3-level* we have two binary terms $Ł_3 = \{101, 111\}$ (decimal 5 and 7), $\#Ł_3 = 2$. And so on. Let's show the first $\#Ł_n$, $n \geq 1$: 1, 1, 2, 3, 6, 10, 20 …

Each level has a *minimum* and a *maximum*. There are no zeros in the binary expansion of maximal terms; for $Ł_n$, $n \geq 1$, these are well-known Mersenne numbers $M_n = \max Ł_n = 2^n - 1$ (the OEIS A000225): 1, 3, 7, 15, 31, 63 … The minimum term of *n*-level, $\min Ł_n$, is defined as the *Dyck Successor*, *DS*, of $\max Ł_{n-1}$ (see [2]): $DS(M_{n-1}) = M_{n-1} + M_m + 1$, $m = \lceil (n-1)/2 \rceil$ ($\lceil \cdot \rceil$ is the ceiling function). We can say that minimal terms are generated by Mersenne numbers. Here are the minima of the initial levels: 1, 3, 5, 11, 19, 39, 71 …

In addition, we know something about some terms within levels. For example, the *n*-level, $n \geq 4$, ends with a *Mersenne triplet* of the form $(M_n - 4, M_n - 2, M_n)$. (Note that in each such triple of odd adjacent numbers all terms are pairwise coprime.) Also, each *n*-level, $n \geq 6$, ends with tree Mersenne triplets (the last nine terms of the *n*-level):

$$(M_n - 20, M_n - 18, M_n - 16), \ (M_n - 12, M_n - 10, M_n - 8), \ (M_n - 4, M_n - 2, M_n).$$

Such nine numbers can be called the *Mersenne nine*. Obviously, in this case, adjacent odd numbers $M_n - 22$, $M_n - 14$, $M_n - 6$, and $M_n + 2$ are not the A036991 terms.

As you know, Dyck paths (and balanced brackets) are counted by the well-known Catalan numbers, the OEIS A000108. However, for Dyck paths of the same length, the corresponding A036991 terms are often at different levels, and vice versa, Dyck paths of different sizes can be encoded by Dyck numbers with the same binary length (this is clearly visible already at the initial levels). In this regard, the problem of enumeration of the A036991/A350346 terms arises, which we will deal with in this paper.

In [3], it is noted that the size of the A036991/A350346 levels is counted by the terms of the OEIS A001405: 1, 1, 2, 3, 6, 10, 20, 35, 70… The following equality is obvious:

(1) $$\#Ł_n = A001405(n-1), \ n \geq 1.$$

Later we will prove the corresponding theorem. Additionally, we consider the concomitant bijection between Dyck numbers and symmetric Dyck paths.

To conclude this section, let's take a look at Figure 1, which clarifies the structure of the OEIS A036991 (we are working with decimal numbers, *tree nodes*). The diagram shows if possible the relationships between the initial levels of the base ternary tree with root 1, as well as the initial nodes of some additional ternary trees. Every ternary tree is infinite, and there are an infinite number of such trees. On the ordinate axis, the level numbers (left column) are set, but on the *x*-axis we have to make a significant offset of the origin at each level, and this allows us to get a compact visual picture.

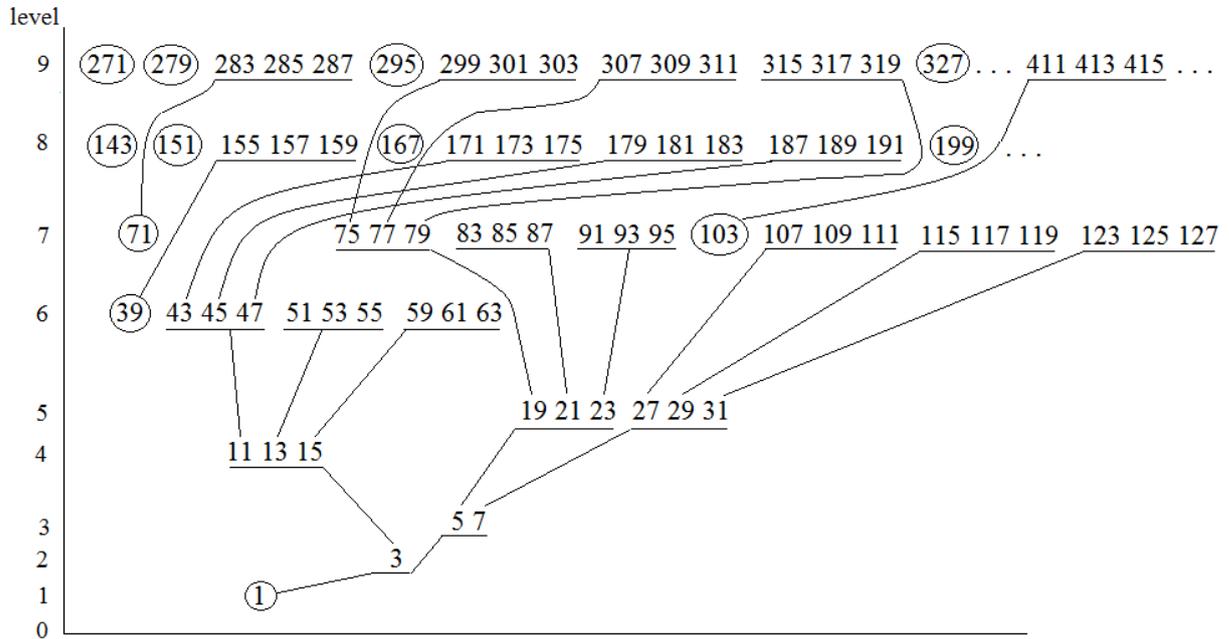

Figure 1. Initial nodes of infinite ternary trees (first nine levels).

We have underlined triples of adjacent odd numbers, *triplets*; the lone terms are shown in a circle, these are the roots of infinite ternary trees. Each node $đ \in ĐN$ (either a tree root or a triplet node) generates a triplet of the form (see [3]): $(4đ – 1, 4đ + 1, 4đ + 3)$. As you can see, a node of *k*-level, $k > 1$, generates a triplet at level $k + 2$. And only at the first level, node 1 generates two-level triplet (3, 5, 7).

The generation of each triplet can be viewed as moving along the ternary tree from the root to infinity. But we can also move backwards, i.e. to the root of the tree: in a triplet, each node $đ$ has the same parent $2 \times \lfloor đ/8 \rfloor + 1$ ($\lfloor \cdot \rfloor$ is the floor function). The root of the ternary tree does not have a parent, and this is a special case.

We can say that starting from 6-level, many tree roots are generated by Mersenne numbers, and we talked about this above. For example, the Mersenne number $M_5 = 31$, the maximum of 5-level, generates a root 39 (31+8), the minimum of 6-level. The Mersenne number $M_7 = 127$, the maximum of 7-level, generates a root 143(+16) at 8-level, and in turn this root generates an additional root 151(+8). At subsequent levels, the number of such roots grows, for example, the 16th Mersenne number 32767 generates a chain of six tree roots: 33023(+256), 33151(+128), 33215(+64), 33247(+32), 33263(+16), 33271(+8). Obviously, the total number of nodes in similar chains is infinite (if we look at all levels).

Of course, the roots of ternary trees are generated inside the levels as well. In Figure 1, these are nodes 103, 167, 199, 295, and 327 (due to lack of space, not all roots of levels 8 and 9 are shown). It is logical to believe that the number of such nodes is also infinite. It is not difficult to draw the obvious conclusion that the number of ternary trees is certainly infinite, and each tree has an infinite number of vertices.

Importantly, ternary trees do not intersect (there are no common nodes), because in each triplet all three nodes have a common (one and only one) parent. Moving through the tree in the opposite direction, we will always get to a single specific root.

## 2  Binary suffixes

We start working with binary suffixes; let us recall the balance rule: in each suffix the number of 0's is no more than 1's. The balance is zero if the suffix has the same number of 0's and 1's. If there are more 1's, the balance is positive; negative balance is invalid (0's greater than 1's). In the case of zero balance, the suffix is a Dyck word, since each Dyck word must have a zero balance.

**2.1. Enumeration.**  In our case, the binary code always ends with 1; otherwise the balance in the smallest one-digit suffix 1, *1-suffix*, is broken. Let us denote the corresponding set $S_1 = \{1\}$. Obviously, the balance of the 1-suffix 1 is 1. We can extend the 1-suffix by adding 0 or 1 first. In both cases, we get valid two-digit suffixes, *2-suffixes*, i.e. $S_2 = \{01, 11\}$. In red, we showed a suffix with zero balance, aka the Dyck word. We can say that we generated the Dyck word 01 by extending the 1-suffix 1.

The next step, the transition from a 2-suffix to a 3-suffix is not obvious: 1 can be added to both 01 and 11, but 0 can only be added to 11. Obviously, 0 cannot be added to Dyck word, otherwise we will get a suffix with a negative balance. As a result, we obtain valid 3-suffixes, $S_3 = \{011, 101, 111\}$. All three new suffixes have an odd length and therefore a positive balance (negative balance is not allowed). Obviously, in the next step, moving from odd suffixes to even suffixes, we can add either 0 or 1 to each odd suffix. As a result, the number of suffixes simply doubles: $S_4 = \{0011, 0101, 0111, 1011, 1101, 1111\}$. Again, two new Dyck words of size 4 are marked in red. This procedure can be continued further.

Table 1 shows suffixes of length $\leq 7$. The first column shows the length of the suffixes, the third column is the cardinality of the corresponding set, $\#S_n$.

| $l$ | Binary suffixes | # |
|---|---|---|
| 0 | $\varepsilon$ | 1 |
| 1 | 1 | 1 |
| 2 | 01, 11 | 2 |
| 3 | 011, 101, 111 | 3 |
| 4 | 0011, 0101, 0111, 1011, 1101, 1111 | 6 |
| 5 | 00111, 01011, 01101, 01111, 10011, 10101, 10111, 11011, 11101, 11111 | 10 |
| 6 | 000111, 001011, 001101, 001111, 010011, 010101, 010111, 011011, 011101, 011111, 100111, 101011, 101101, 101111, 110011, 110101, 110111, 111011, 111101, 111111 | 20 |
| 7 | 0001111, 0010111, 0011011, 0011101, 0011111, 0100111, 0101011, 0101101, 0101111, 0110011, 0110101, 0110111, 0111011, 0111101, 0111111, 1000111, 1001011, 1001101, 1001111, 1010011, 1010101, 1010111, 1011011, 1011101, 1011111, 1100111, 1101011, 1101101, 1101111, 1110011, 1110101, 1110111, 1111011, 1111101, 1111111 | 35 |

Table 1. Binary suffixes of length $\leq 7$.

Empty words are often used, we also added the empty suffix $\varepsilon$ (of size 0), $\#S_0 = 1$. Since the empty word is Dyck word, we also marked it in red. Now let us explain the presence of a single 1-suffix. We can add 1 to the empty word $1 \cdot \varepsilon = 1$, but adding zero gives an invalid 1-suffix $0 \cdot \varepsilon = 0$, so $\#S_1 = 1$ ($\cdot$ is a concatenation operation).

As you can see, Dyck words appear among the even-length suffixes, and their number gradually grows: 1, 1, 2, 5 (obviously, at the next even ranges we get 14, 42, 132, and so on). And this is natural, since when expanding binary suffixes, we generate Dyck words that are counted by these Catalan numbers (the OEIS A000108). Recall that the number of Dyck words with length of $2n$ is the $n$th Catalan number $Cat(n) = (2n)!/(n!(n+1)!)$. We can make the following obvious statement.

**Proposition 2**. Let $S_n$ be the set of all valid suffixes of length $n$, $n \geq 0$. Then:

*(i)*  In the case of odd $n$, $\#S_{n+1} = 2 \times \#S_n$.
*(ii)* In the case of even $n$ (including 0), $\#S_{n+1} = 2 \times \#S_n - Cat(n/2)$.

In Table 1, the numbers in the right column match the A001405 terms, and this is simply explained: if we add 1 in front of each valid $k$-suffix, we get all the A350346 terms of $(k+1)$-level (see equality (1)). The following equality is true:

(2) $$\#S_n = A001405(n), \ n \geq 0.$$

Let us formulate an obvious lemma.

**Lemma 3**. The number of Dyck numbers of binary length $n$, $n > 0$, is equal to the number of admissible binary suffixes of length $n-1$, i.e.

(3) $$\#Ł_n = \#S_{n-1}, \ n \geq 1.$$

**2.2. Generating function.** In the right column of Table 1, we actually have two sequences. One sequence is formed at even levels (suffixes of even length): 1, 2, 6, 20, 70, 262 ..., and these are the central binomial coefficients $(2n!)/(n!)^2$, $n \geq 0$ (see the OEIS A000984). The generating function of this sequence is well known [5, 6]:

$$G(x) = 1/(1-4x)^{1/2} = 1 + 2x + 6x^2 + 20x^3 + 70x^4 + 252x^5 + ...$$

Another sequence is formed at odd levels (odd suffixes): 1, 3, 10, 35, 126... The generating function of this sequence is usually obtained as follows (see the OEIS A001700):

$$H(x) = (G(x) - 1)/2x = (1/(1-4x)^{1/2} - 1)/2x = 1 + 3x + 10x^2 + 35x^3 + 126x^4 + ...$$

It remains to nest both sequences one into the other, distributing the coefficients over the corresponding powers of $x$.

$$G(x^2) = 1/(1-4x^2)^{1/2} = 1 + 2x^2 + 6x^4 + 20x^6 + 70x^8 + 252x^{10} + ...$$
$$xH(x^2) = x(G(x^2) - 1)/(2x^2)) = (1/(1-4x^2)^{1/2} - 1)/(2x)$$
$$= x + 3x^3 + 10x^5 + 35x^7 + 126x^9 + ...$$
$$G(x^2) + xH(x^2) = 1/(1-4x^2)^{1/2} + (1/(1-4x^2)^{1/2} - 1)/(2x)$$
$$= 1 + x + 2x^2 + 3x^3 + 6x^4 + 10x^5 + 20x^6 + 35x^7 + 70x^8 \ ...$$

# 3    Symmetric Dyck numbers, the bijection $\mathcal{B}$

According to equalities (2) and (3), the Dyck numbers and their suffixes are counted by the A001405 terms. In addition, we note that in the OEIS A001405 there are comments on the listing of symmetrical Dyck paths (Matt Watson, 2012) as well as symmetric balanced parentheses (Joerg Arndt, 2011). In this section, we will construct the corresponding bijection on the Dyck numbers.

Symmetry is well seen both on Dyck paths and on balanced parentheses, but the brackets are more compact (there are no bulky graphs), so we will work with bracket sets. The concept of symmetry is quite obvious. For example, the brackets (()()) is symmetric and the brackets (())() are not. In this regard, we define the corresponding Dyck numbers.

**Definition 4.** Let us call a Dyck number *symmetrical* if it encodes a symmetric Dyck path. Otherwise the Dyck number is *asymmetrical*.

From the example above, the first balanced brackets are encoded with the integer $11 = 001011_2$, which is a symmetrical Dyck number, while $13 = 001101_2$ is not. The first binary code is symmetrical and the second one is asymmetric. In A036991 there are significantly more asymmetric terms at high levels (20 and above), but in the initial ranges the symmetric Dyck numbers predominate. Here are the first symmetric terms of the A036991 (currently there is no such sequence in the OEIS):

(4)    1, 3, 5, 7, 11, 15, 21, 23, 31, 43, 47, 51, 63, 77, 85, 87, 95, 103, 127, 155, 171 …

As you can see, all Mersenne numbers are symmetrical. With increasing ranges, the relative number of symmetric terms tends to zero. For example, in the 18th level, out of 24310 terms, only 190 are symmetrical (less than one percent). In Appendix A we give symmetric A036991 terms for the first 16 levels.

Below we will often use the asymmetric A036991 terms as well. Here are the first asymmetric terms (and this sequence is not yet in the OEIS):

(5)    13, 19, 27, 29, 39, 45, 53, 55, 59, 61, 71, 75, 79, 83, 91, 93, 107, 109, 111, 115 …

Asymmetric terms are very important in this paper. Below we consider unary trees, the roots of trees are asymmetric terms. Appendix B shows the asymmetric terms for the first 11 levels.

**Example 5.** The OEIS A001405 contains important comments and accompanying examples. Emeric Deutsch (2005) notes that the *n*th term is equal to the number of left factors of Dyck paths, consisting of *n* steps. Prefixes are given that correspond to a(4) = 6: UUDD, UDUD, UUUD, UUDU, UDUU, and UUUU, where U = (1, 1) and D = (1, -1). Lee A. Newberg (2010) gives similar bracket prefixes of length 4: (()), ()(), (((, (()(, ()((, and ((((. We are working with binary suffixes, so in our case let us show the corresponding 4-suffixes from Table 1: 0011, 0101, 0111, 1011, 1101, and 1111. And if we assign 1 to the beginning of each binary suffix (equivalent to add $2^4 = 16$), we get the following six Dyck numbers: 19, 21, 23, 27, 29, and 31.

Obviously, to move from an arbitrary path prefix to a symmetric path, you just need to invert the prefix into a suffix of the same length and combine (concatenate) both frag-ments. We show this procedure on the mentioned six bracket prefixes (the inverted fragments are shown in red):

(()) → (())(()), ()() → ()()()(), ((() → ((()())), ()(( → ()(()()), )(( → ()(())(), (((( → ((((())))

For binary suffixes, let's show the entire chain of transformations, starting from the original A036991 term (decimal number is in brackets) and ending with the resulting symmetric Dyck number (again, the inverted fragments are shown in red).

(19) 10011 → 0011 → 0011·0011 = 00110011 (51),
(21) 10101 → 0101 → 0101·0101 = 01010101 (85),
(23) 10111 → 0111 → 0001·0111 = 00010111 (23),
(27) 11011 → 1011 → 0010·1011 = 00101011 (43),
(29) 11101 → 1101 → 0100·1101 = 01001101 (77),
(31) 11111 → 1111 → 0000·1111 = 00001111 (15).

To conclude Example 5, let us fix that first (*i*) in the binary we remove the leading bit, the 1, then (*ii*) we invert the resulting suffix and obtain a prefix of the same length, then (*iii*) the prefix is concatenated with the original suffix, and finally, (*iv*) after removing the leading zeros, we obtain the symmetric A036991 term.

The four-step procedure considered in Example 5 can be seen as an *algorithm* for obtaining the symmetric term A036991. In the six chains shown, we obtained a kind of bijection (one-to-one correspondence function) between arbitrary Dyck numbers and symmetric Dyck numbers. All terms of the OEIS A036991 (except 0) start in binary with 1; by removing this bit, we actually dump the ballast, while reducing the resulting integers. The uniqueness of the A036991 term with binary length *n*, *n* > 0, automatically guarantees the uniqueness of its binary (*n*-1)-suffix. Obviously, the reverse of such a suffix and its subsequent concatenation with the resulting (*n*-1)-prefix does not break the bijection.

Let us denote the resulting bijection as $\mathcal{B}: \text{Đ}\mathbb{N}\setminus\{0\} \rightarrow s\text{Đ}\mathbb{N}$, where sĐN is the set of symmetric Dyck numbers (the numerical codes of symmetric Dyck paths and symmetric balanced brackets). Obviously, bijection $\mathcal{B}$ transforms an arbitrary *n*-level A036991 term into a term that encodes a symmetric Dyck path of semilength *n*-1. Table 2 shows in detail the bijection $\mathcal{B}$ in the initial levels.

| *đ* | Binary codes and brackets | Suffixes | Symmetry | $\mathcal{B}(đ)$ |
|---|---|---|---|---|
| **1** | 1 = 01 → () | ε | | 0 |
| **3** | 11 = 0011 → (()) | 1 → ) | () → 01 | **1** |
| **5** | 101 = 0101 → ()() | 01 → () | ()() → 0101 | 5 |
| **7** | 111 = 000111 → ((())) | 11 → )) | (()) → 0011 | **3** |
| 11 | 1011 = 001011 → (()()) | 011 → ()) | (()()) → 001011 | 11 |
| 13 | 1101 = 001101 → (())() | 101 → )() | ()()() → 010101 | 21 |
| **15** | 1111 = 00001111 → (((()))) | 111 → ))) | ((())) → 000111 | **7** |
| 19 | 10011 = 010011 → ()(()) | 0011 → (()) | (())(()) → 00110011 | 51 |
| 21 | 10101 = 010101 → ()()() | 0101 → ()() | ()()()() → 01010101 | 85 |
| 23 | 10111 = 00010111 → ((()())) | 0111 → ()) | ((()())) → 00010111 | 23 |
| 27 | 11011 = 00011011 → ((())()) | 1011 → )()) | (()()()) → 00101011 | 43 |
| 29 | 11101 = 00011101 → ((()))() | 1101 → ))() | ()(())() → 01001101 | 77 |
| **31** | 11111 = 0000011111 → (((()))) | 1111 → )))) | ((((())))) → 00001111 | **15** |
| 39 | 100111 = 00100111 → (())(()) | 00111 → (()) | ((())(())) → 0001100111 | 103 |
| 43 | 101011 = 00101011 → (()()()) | 01011 → ()()) | (()()()) → 0010101011 | 171 |
| 45 | 101101 = 00101101 → (())()() | 01101 → ())() | ()(())()() → 0100101101 | 301 |
| 47 | 101111 = 0000101111 → (((()()))) | 01111 → ())))) | (((()()))) → 0000101111 | 47 |
| 51 | 110011 = 00110011 → (())(()) | 10011 → )(()) | (())()(()) → 0011010011 | 211 |
| 53 | 110101 = 00110101 → (())()() | 10101 → )()() | ()()()()() → 0101010101 | 341 |
| 55 | 110111 = 0000110111 → ((()())()) | 10111 → )())) | (((()())))→ 0001010111 | 87 |
| 59 | 111011 = 0000111011 → ((((())())) | 11011 → ))()) | (()(()))()) → 0010011011 | 155 |
| 61 | 111101 = 0000111101 → (((()))() | 11101 → )))() | ()((())))() → 0100011101 | 285 |
| **63** | 111111 = 000000111111 → (((((()))))) | 11111 → ))))) | (((((())))))→ 0000011111 | **31** |

Table 2. The bijection $\mathcal{B}$ for the first six levels.

It is logical to consider 0 a symmetrical Dyck number; 0 encodes an empty Dyck word $\varepsilon$, and the equality $\mathcal{B}(1) = \mathcal{B}(M_1) = M_0 = 0$ is obvious. In general, for Mersenne numbers, we have $\mathcal{B}(M_n) = M_{n-1}$, $n > 0$.

The following terms are interesting (OEIS sequence [A052940](A052940) without initial term)

(6)        5, 11, 23, 47, 95, 191, 383, 767, 1535, 3071, 6143, 12287, 24575, 49151 ...

with equalities $a_n = 3 \times 2^n - 1$ and $a_{n+1} = 2a_n + 1$ (the recurrence relation is similar to Mersenne numbers). The bijection $\mathcal{B}$ is the *identity* on the terms of (6). These terms (shown in blue in Table 2) can be called *self-bijective*, since $\mathcal{B}(a_n) = \mathcal{B}^{-1}(a_n) = a_n$.

Bijection $\mathcal{B}$ can be attributed to trivial bijective functions, since the set of asymmetric Dyck numbers ĐN\sĐN is countable and $|ĐN\setminus\{0\}| = |sĐN|$. For a symmetrical term $đ$ it is not difficult to obtain the inverse term $\mathcal{B}^{-1}(đ)$. Let's formulate a common definition.

**Definition 6.** For an arbitrary natural number $n$, the *binary weight*, bw($n$), is the sum of the digits in the binary expansion.

Any A036991 term $đ > 1$ begins and ends with 1 in binary expansion, so bw($đ$) $\geq 2$, $đ > 1$. Obviously, for Mersenne numbers, bw($M_n$) = $n$.

**Proposition 7.** Let $đ$ be a symmetric A036991 term, then

(7) $$\mathcal{B}^{-1}(đ) = (đ \bmod 2^w) + 2^w, \text{ where } w = \text{bw}(đ).$$

*Proof.* Obviously, $w$ is equal to the semilength of symmetric Dyck path, which is encoded by $đ$. Such a path has $w$ up-steps and $w$ down-steps. Therefore, the operation of taking modulo $2^w$ allows us to extract the required binary suffix of length $w$ (in binary, this also includes possible leading zeros). The resulting binary suffix is equal to half of the symmetric path. The addition operation with $2^w$ adds 1 at the beginning of the binary suffix. The result shows that $\mathcal{B}$ is indeed a bijection. For example, if $đ = 0$, then bw(0) = 0 and

$$\mathcal{B}^{-1}(0) = (0 \bmod 2^0) + 2^0 = 0 + 1 = 1. \qquad \square$$

According to formula (7) for the Mersenne numbers, we get

$$\mathcal{B}^{-1}(M_n) = (M_n \bmod 2^n) + 2^n = M_n + 2^n = M_{n+1}.$$

Thus, for Mersenne numbers the inverse bijective function, $\mathcal{B}^{-1}$, gives us the A000225: $M_{n+1} = \mathcal{B}^{-1}(M_n)$, $n \geq 0$. But the original function $\mathcal{B}$ for a given A000225 term generates a finite set of Mersenne numbers ending in $M_0 = 0$.

# 4 Unary trees of the bijection $\mathcal{B}$

In the case of Mersenne numbers, the functions $\mathcal{B}$ and $\mathcal{B}^{-1}$ deal only with symmetric terms. The situation is similar with the symmetric terms of the infinite sequence (6), on which the bijection works as an identity function (or identity transformation).

Let us look at how the bijection $\mathcal{B}$ works on the remaining A036991 terms, that is, on the set $ĐN_0 = (ĐN \setminus A000225) \setminus A052940$.

In Table 3, the left column shows the initial asymmetric A036991 terms (root of unary trees), and the rows show the chains of symmetrical terms, which are obtained by the bijection $\mathcal{B}$ gradually step-by-step, starting with the asymmetric terms from the left column.

| Root | Chain, sequential terms of the bijection $\mathcal{B}$ |
|------|---------------------------------------------------------|
| 13   | 21, 85, 1365, 349525, 22906492245, 98382635059784275285, … |
| 19   | 51, 211, 3411, 873811, 57266230611, 245956587649460688211, … |
| 27   | 43, 171, 2731, 699051, 45812984491, 196765270119568550571, … |
| 29   | 77, 1229, 314573, 20615843021, 88544371553805847757, … |
| 39   | 103, 423, 6823, 1747623, 114532461223, 491913175298921376423, … |
| 45   | 301, 19245, 78826285, 1322485604862765, 3722466048289245067887144333325, … |
| 53   | 341, 21845, 89478485, 1501199875790165, 4225502000760764671655677335125, … |
| 55   | 87, 343, 5463, 1398103, 91625968983, 393530540239137101143, … |
| 59   | 155, 2459, 629147, 41231686043, 177088743107611695515, … |
| 61   | 285, 18205, 74565405, 1250999896491805, 35212516673006372263797311260, … |
| 71   | 455, 7367, 1887431, 123695058119, 531266229322835086535, … |
| 75   | 715, 45771, 187478731, 3145371168322251, 88534327634987450263261811169, … |
| 79   | 207, 847, 13647, 3495247, 229064922447, 983826350597842752847, … |
| 83   | 851, 54611, 223696211, 3752999689475411, 1056375500190191167913919337811, … |
| 91   | 603, 38491, 157652571, 2644971209725531, 7444932096578490135774288666651, … |
| 93   | 1117, 285789, 18729426013, 8044227219594126038l, … |

Таблица 3. Successive symmetric terms of the bijection $\mathcal{B}$ with tree roots $< 100$.

For an asymmetric term $đ \in ĐN_0$ we have shown the beginning of an infinite chain of symmetric terms of the form $\mathcal{B}(đ)$, $\mathcal{B}(\mathcal{B}(đ))$, $\mathcal{B}(\mathcal{B}(\mathcal{B}(đ)))$, $\mathcal{B}(\mathcal{B}(\mathcal{B}(\mathcal{B}(đ))))$, and so on. For example, for the first row we have $\mathcal{B}(13) = 21$, $\mathcal{B}(21) = 85$, $\mathcal{B}(85) = 1365$… We can say that we got an infinite unary tree, which is generated by the bijection $\mathcal{B}$; the root of this tree is the asymmetric term 13. Let us call such a tree *bijective* or *bijection tree*, $\mathcal{B}$-*tree*.

If in each row of Table 3 we follow the chain in the opposite direction using the function $\mathcal{B}^{-1}$, we will return to a certain asymmetric term, the root of the corresponding bijection tree. Obviously, such $\mathcal{B}$-trees are infinite, do not intersect, and there are an infinite number of bijection trees. Let us formulate an appropriate statement.

**Proposition 8**. The bijection $\mathcal{B}$ splits the set $ĐN_0$ into an infinite number of disjoint bijection trees. The roots of the trees are all asymmetric A036991 terms, the other nodes are symmetric terms of $ĐN_0$ (not Mersenne numbers and not A052940 terms).

Clearly, the set $ĐN_0$ is a forest, and in this forest, at first glance, strange things happen.

**Contradiction on a countable set.** The question may arise, how many of Dyck paths are symmetrical and how many are asymmetrical? And which ones are more? For example, in the first 20 levels of the OEIS A036991, the proportion of symmetric terms is about one percent of all terms. And this ratio decreases significantly with each successive level; already in 28 levels, the share of symmetric terms is 0.033% (the relative number of symmetric terms confidently tends to zero).

But on the other hand, each bijection tree has only one asymmetric term, and that is the root of the tree. The remaining nodes of the $\mathcal{B}$-tree are symmetrical, and there are an infinite number of them (recall that the bijection $\mathcal{B}$ does not generate asymmetric terms). The number of bijection trees is also infinite. From this we draw the opposite conclusion: the number of asymmetric A036991 terms is negligible. (Recall that outside bijection trees there are also symmetric Mersenne numbers and symmetric A052940 terms.)

Probably, such a contradiction often arises on countable sets. In the countable set of natural numbers, there is an infinite subset of composite numbers and an infinite subset of primes (Euclid's fundamental theorem). Both subsets do not intersect, there are no common elements. As we know, there are significantly fewer primes (on a finite segment of natural numbers). Both subsets are countable, and therefore there must be a bijection between them. But so far no such bijection has been discovered, and we think no one is looking for it.

In our case, it is also problematic to construct a bijection between countable subsets of the set $ĐN_0$ (between asymmetric terms and symmetric terms). The bijection trees have no common nodes, so the replenishment of chains (see Table 3) is achieved by selecting symmetric terms from very distant levels. Let us try to clarify the situation with an example.

**Example 9.** In the bijection tree with the asymmetric root 45 (the sixth level in the OEIS A036991) the fifth symmetric term 3722466048289245067887144333325 is 30 digits long, and the 8th term already has 232 decimal digits. This term is

7294370719381108647684094104522234227873623191532415823764434972312874317383612037364127127194910436084831242640488851889314090365999976285503982592717761297150029639527256463560124642138020856616314044065619166270 2607661408405605165.

The 8th term is at level 771, which corresponds to its length in binary expansion. In bijection tree chains, the length of each new term almost doubles, so the binary length of the next 9th term of the chain with root 45 will exceed 1500 digits (and this is already a giant level in A036991).

Each $\mathcal{B}$-tree is infinite, the number of such trees is infinite, so the OEIS A036991, which includes this entire forest, is also infinite. And it is logical that the replenishment of tree chains comes from the depths of A036991. Obviously, any two terms from the same level cannot fall into the same chain. But each term from ÐN₀ is guaranteed to be in a chain and only one chain of a certain $\mathcal{B}$-tree.

It is a good idea to check, at least at the levels available to us, whether all the initial symmetric terms of ÐN₀ are disassembled by chains of bijection trees. Let us look at the following example of how the chains are filled at the first levels (see Table 3).

**Example 10.** In the first 9 levels of the OEIS A036991 there are 21 symmetric Dyck numbers (see Appendix A), which are included in the chains of bijection trees:  21, 43, 51, 77, 85, 87, 103, 155, 171, 175, 207, 211,  285, 301, 311, 341, 343, 351, 415, 423, 455.  17 numbers (shown in red) are selected by asymmetric terms directly in Table 3. Here are the initial parts of such chains (we have underlined the roots of trees):

13↦21↦85, 19↦51↦211, 27↦43↦171, 29↦77, 39↦103↦423, 45↦301, 53↦341, 55↦87↦343, 59↦155, 61↦285, 71↦455, 79↦207.

For the other four Dyck numbers, tree roots can be obtained manually (but this is quite painstaking work):  111↦175, 119↦311, 159↦415, 223↦351.

We did this programmatically; the corresponding programs are shown in the supplementary applications. The first Python program (Appendix C) allows you to obtain chains for bijection trees in a certain range of roots. The second Python program (Appendix D) gets a symmetric Dyck number from the operator and prints in reverse order the preceding terms of the chain, ending with the tree root. The second program performs the inverse bijective function  $\mathcal{B}^{-1}$. That's exactly how we got the last four roots.

**Conclusion**. The considered bijection $\mathcal{B}$ and unary $\mathcal{B}$-trees allow us to look into the depths of the OEIS A036991 and accordingly push us to study the structure of this sequence more thoroughly and in more depth. Unfortunately, the initial levels provide little information, due to the complex A036991 structure.

Acknowledgements. The author would like to thank Olena G. Kachko (Kharkiv National University of Radio Electronics) for helpful discussions.

# References

[1]  R. Stanley. *Catalan numbers*. Cambridge University Press, Cambridge, 2015.

[2]  G. Eremin. *Dyck numbers, I. Successor function*, 2022. arXiv:2210.00744

[3]  G. Eremin. *Dyck numbers, II. Triplets and rooted trees in OEIS A036991*, 2022. arXiv:2211.01135

[4]  N. J. A. Sloane, The On-Line Encyclopedia of Integer Sequences, https://oeis.org.

[5]  Li-Hua Deng, Eva Y. P. Deng and Louis W. Shapiro. *The Riordan Group and Symmetric Lattice Paths*, 2009. arXiv:0906.1844

[6]  Paul Barry. *The Central Coefficients of a Family of Pascal-like Triangles and Colored Lattice Paths*, J. Int. Seq., Vol. 22 (2019), Article 19.1.3. https://www.emis.de/journals/JIS/VOL22/Barry1/barry411.pdf

Concerned with sequences A000108, A000225, A000984, A001405, A001700, A036991, A052940, and A350346.

# A  The symmetric A036991 terms without Mersenne numbers and A052940 terms (total 333)

| Level | Symmetric terms | # |
|---|---|---|
| 5 | 21 | 1 |
| 6 | 43 51 | 2 |
| 7 | 77 85 87 103 | 4 |
| 8 | 155 171 175 207 211 | 5 |
| 9 | 285 301 311 341 343 351 415 423 455 | 9 |
| 10 | 571 603 623 683 687 703 715 819 831 847 851 911 | 12 |
| 11 | 1085 1117 1143 1197 1207 1229 1247 1333 1365 1367 1375 1407 1431 1639 1663 1695 1703 1823 1863 | 19 |
| 12 | 2171 2235 2287 2395 2415 2459 2495 2667 2731 2735 2751 2815 2863 2891 3187 3251 3279 3327 3391 3407 3411 3647 3727 3855 | 24 |
| 13 | 4221 4285 4343 4445 4471 4509 4575 4717 4781 4791 4831 4919 4941 4991 5237 5301 5335 5461 5463 5471 5503 5631 5727 5783 5911 6375 6503 6559 6655 6783 6815 6823 6951 7295 7367 7455 7495 7711 | 38 |
| 14 | 8443 8571 8687 8891 8943 9019 9151 9435 9563 9583 9663 9839 9883 9983 10011 10475 10603 10671 10923 10927 10943 11007 11051 11263 11455 11467 11567 11595 11823 12531 12659 12751 12979 13007 13107 13119 13311 13523 13567 13631 13647 13651 13903 14591 14735 14911 14991 15423 15631 | 49 |
| 15 | 16637 16765 16887 17085 17143 17213 17375 17629 17757 17783 17887 18039 18077 18205 18303 18669 18797 18871 19117 19127 19167 19245 19327 19661 19679 19767 19789 19967 20023 20725 20853 20951 21173 21207 21301 21343 21717 21845 21847 21855 21887 22015 22103 22527 22911 22935 23135 23191 23647 23831 25063 25319 25503 25959 26015 26215 26239 26623 27047 27135 27263 27295 27303 27807 27943 29127 29183 29383 29471 29823 29983 30023 30847 31263 31775 | 75 |
| 16 | 33275 33531 33775 34171 34287 34427 34751 35259 35515 35567 35775 36079 36155 36411 36607 37339 37595 37743 38235 38255 38335 38491 38655 39323 39359 39535 39579 39935 40047 40219 41451 41707 41903 42347 42415 42603 42687 43435 43691 43695 43711 43775 44031 44207 44331 45055 45515 45771 45823 45871 46271 46383 46411 47295 47663 48175 49651 49907 50127 50547 50639 50803 51007 51635 51891 51919 52031 52431 52479 52531 53247 53715 53971 54095 54271 54527 54591 54607 54611 55615 55887 56399 58255 58367 58767 58943 59647 59967 60047 60559 61695 62223 62527 62735 63551 | 95 |

# B  The asymmetric A036991 terms (total 454)

| Level | Asymmetric A036991 terms | # |
|---|---|---|
| 4 | 13 | 1 |
| 5 | 19, 27, 29 | 3 |
| 6 | 39, 45, 53, 55, 59, 61 | 6 |
| 7 | 71, 75, 79, 83, 91, 93, 107, 109, 111, 115, 117, 119, 123, 125 | 14 |
| 8 | 143, 151, 157, 159, 167, 173, 179, 181, 183, 187, 189, 199, 203, 205, 213, 215, 219, 221, 223, 231, 235, 237, 239, 243, 245, 247, 251, 253 | 28 |
| 9 | 271, 279, 283, 287, 295, 299, 303, 307, 309, 315, 317, 319, 327, 331, 333, 335, 339, 347, 349, 359, 363, 365, 367, 371, 373, 375, 379, 381, 399, 407, 411, 413, 427, 429, 431, 435, 437, 439, 443, 445, 447, 459, 461, 463, 467, 469, 471, 475, 477, 479, 487, 491, 493, 495, 499, 501, 503, 507, 509 | 59 |
| 10 | 543, 559, 567, 573, 575, 591, 599, 605, 607, 615, 619, 621, 627, 629, 631, 635, 637, 639, 655, 663, 667, 669, 671, 679, 685, 691, 693, 695, 699, 701, 711, 717, 719, 723, 725, 727, 731, 733, 735, 743, 747, 749, 751, 755, 757, 759, 763, 765, 783, 791, 795, 797, 799, 807, 811, 813, 815, 821, 823, 827, 829, 839, 843, 845, 853, 855, 859, 861, 863, 871, 875, 877, 879, 883, 885, 887, 891, 893, 895, 919, 923, 925, 927, 935, 939, 941, 943, 947, 949, 951, 955, 957, 959, 967, 971, 973, 975, 979, 981, 983, 987, 989, 991, 999, 1003, 1005, 1007, 1011, 1013, 1015, 1019, 1021 | 112 |
| 11 | 1055, 1071, 1079, 1083, 1087, 1103, 1111, 1115, 1119, 1127, 1131, 1133, 1135, 1139, 1141, 1147, 1149, 1151, 1167, 1175, 1179, 1181, 1183, 1191, 1195, 1199, 1203, 1205, 1211, 1213, 1215, 1223, 1227, 1231, 1235, 1237, 1239, 1243, 1245, 1255, 1259, 1261, 1263, 1267, 1269, 1271, 1275, 1277, 1279, 1295, 1303, 1307, 1309, 1311, 1319, 1323, 1325, 1327, 1331, 1335, 1339, 1341, 1343, 1351, 1355, 1357, 1359, 1363, 1371, 1373, 1383, 1387, 1389, 1391, 1395, 1397, 1399, 1403, 1405, 1423, 1435, 1437, 1439, 1447, 1451, 1453, 1455, 1459, 1461, 1463, 1467, 1469, 1471, 1479, 1483, 1485, 1487, 1491, 1493, 1495, 1499, 1501, 1503, 1511, 1515, 1517, 1519, 1523, 1525, 1527, 1531, 1533, 1567, 1583, 1591, 1595, 1597, 1599, 1615, 1623, 1627, 1629, 1631, 1643, 1645, 1647, 1651, 1653, 1655, 1659, 1661, 1679, 1687, 1691, 1693, 1707, 1709, 1711, 1715, 1717, 1719, 1723, 1725, 1727, 1735, 1739, 1741, 1743, 1747, 1749, 1751, 1755, 1757, 1759, 1767, 1771, 1773, 1775, 1779, 1781, 1783, 1787, 1789, 1791, 1807, 1815, 1819, 1821, 1831, 1835, 1837, 1839, 1843, 1845, 1847, 1851, 1853, 1855, 1867, 1869, 1871, 1875, 1877, 1879, 1883, 1885, 1887, 1895, 1899, 1901, 1903, 1907, 1909, 1911, 1915, 1917, 1919, 1935, 1943, 1947, 1949, 1951, 1959, 1963, 1965, 1967, 1971, 1973, 1975, 1979, 1981, 1983, 1991, 1995, 1997, 1999, 2003, 2005, 2007, 2011, 2013, 2015, 2023, 2027, 2029, 2031, 2035, 2037, 2039, 2043, 2045 | 231 |

# C  Python program for building unary trees in ĐN₀

```python
# Online Python - IDE, Editor, Compiler, Interpreter
def ok(n):     # Dyck number selection
    if n < 1: return 0
    count = {'0': 0, '1': 0}
    for bit in bin(n)[:1:-1]:
        count[bit] += 1
        if count['0'] > count['1']: return 0
    return count['1']-count['0'] # the number of leading 0's
def ussym(nb):    # transformation of the Dyck number
    dbl = []       # into a symmetric number
    for byte in nb[::-1]:
        if byte=='1': dbl.append('0')
        else: dbl.append('1')
    dbl.pop(); return ''.join(dbl) + nb[1:]

def is_symdyck(nb):  # test of the Dyck number for symmetry
    i= 0; j= len(nb) - 1
    while i < j:
        if nb[i] == nb[j]: return 0 #  not symmetrical
        i+= 1; j-= 1;
    return 1
def dn_chain(sd, msymm):   # creation of a tree chain
    msdn= [sd]
    for k in range(3):    # in chain sd and 3 new terms
        bnydn = bin(sd)[2:]; bsm = ussym(bnydn);
        sd = int(bsm, 2); msdn.append(sd);
        if sd in msymm: msymm.remove(sd)
    return msdn

masym= []; msymm= []
for n in range(3, 90000, 2):
    zers= ok(n)
    if zers==0: continue
    bindn = bin(n)[2:]; d00 = '0'*zers + bindn
    if is_symdyck(d00): msymm.append(n)
    else:
        if len(masym) < 7000: masym.append(n) # try 5000
print('Symmetry', len(msymm), msymm, '\n'); j = 0
while j < len(masym):
    k = masym[j]; bindn = bin(k)[2:]
    bsym = ussym(bindn); sdn = int(bsym, 2)
    if sdn in msymm: msymm.remove(sdn)
    print(k, '-->', dn_chain(sdn, msymm) ); j += 1

print('Not selected', msymm); # Mersenne numbers and A052940
```

# D     Python program for the inverse of the bijection $\mathcal{B}$

```python
# Online Python - IDE, Editor, Compiler, Interpreter
def ok(n):     # Dyck number selection
    if n < 1: return 0
    count = {'0': 0, '1': 0}
    for bit in bin(n)[:1:-1]:
        count[bit] += 1
        if count['0'] > count['1']: return 0
    return count['1']-count['0'] # the number of leading 0's
def ussym(nb): # transformation of Dyck number into symmetric
    dbl = []
    for byte in nb[::-1]:
       if byte=='1': dbl.append('0')
       else: dbl.append('1')
    dbl.pop(); return ''.join(dbl) + nb[1:]
def is_symdyck(nb):
    i= 0; j= len(nb) - 1
    while i < j:
        if nb[i] == nb[j]: return 0
        i+= 1; j-= 1;
    return 1
def inverse_bij(dn):
    hlf = len(dn)//2; return '1' + dn[hlf:]
mers= [3, 7, 15, 31, 63, 127, 255, 1023, 2047, 4095, 8191, 16383]
a52940= [11, 23, 47, 95, 191, 383, 767, 1535, 3071, 6143, 12287, 24575, 49151, 98303, 196607, 393215, 786431, 1572863, 3145727]
n = int(input('Please symmetric A036991 term: '))
while 1:
    if ok(n)==0: print('It is not A036991 term'); break
    if n in mers: print(n, '-->', (2*n+1)); break
    if n in a52940: print(n, '-->', n); break
    chain= [n]; nn= n
    while 1:
        bindn = bin(nn)[2:];
        bdn0= '0'*ok(nn) + bindn # with leading 0's
        if is_symdyck(bdn0)==0: break    # asymmetric DN
        nbij = inverse_bij(bdn0); nij= int(nbij, 2)
        nn = nij; chain.append(nn)
    if len(chain)>1: print(chain); break
    print(n, 'is a tree root, asymmetric term'); break
```